\documentclass[12pt]{amsart}
\usepackage{amsmath,amssymb,amsthm,amscd}
\newtheorem{thm}{Theorem}
\newtheorem{lem}[thm]{Lemma}

\newtheorem*{rem}{Remark}
\newtheorem*{ack}{Acknowledgements}

\def\C{\mathbb{C}}
\def\Z{\mathbb{Z}}
\def\Q{\mathbb{Q}}
\def\R{\mathbb{R}}
\def\F{\mathbb{F}}

\def\N{\mathbb{N}}

\newcommand{\legen}[2]{\genfrac{(}{)}{}{}{#1}{#2}}

\def\ord{{\rm ord}}

\def\SL{{\rm SL}}
\newcommand{\tors}{\operatorname{tors}}
\newcommand{\rank}{\operatorname{rank}}

\usepackage[OT2,T1]{fontenc}
\DeclareSymbolFont{cyrletters}{OT2}{wncyr}{m}{n}
\DeclareMathSymbol{\Sha}{\mathalpha}{cyrletters}{"58}

\begin{document}

\author{Marvin Jones}
\address{Department of Mathematics, University of South Carolina, Columbia,
SC 29208}
\email{jonesmc8@gmail.com}

\author{Jeremy Rouse}
\address{Department of Mathematics, Wake Forest University, Winston-Salem,
NC 27109}
\email{rouseja@wfu.edu}

\title{Solutions of the cubic Fermat equation in quadratic fields}

\begin{abstract}
We give necessary and sufficient conditions on a squarefree integer $d$
for there to be non-trivial solutions to $x^{3} + y^{3} = z^{3}$
in $\Q(\sqrt{d})$, conditional on the Birch and
Swinnerton-Dyer conjecture. These conditions are similar to those obtained by
J. Tunnell in his solution to the congruent number problem.
\end{abstract}


\maketitle

\section{Introduction and Statement of Results}

The enigmatic claim of Fermat that the equation
\[ 
  x^{n} + y^{n} = z^{n}
\]
has only the trivial solutions (those with at least one of $x$, $y$ and $z$ zero)
in integers when $n \geq 3$ has to a large extent shaped 
the development of number theory over the course of the last three hundred 
years. These developments culminated in the theory used by
Andrew Wiles in \cite{Wiles} to finally justify Fermat's claim.

In light of Fermat's claim and Wiles's proof, it is natural to ask the
following question:  for which fields $K$ does the equation $x^{n} +
y^{n} = z^{n}$ have a non-trivial solution in $K$? Two notable results
on this question are the following. In \cite{JM}, it is shown that the
equation $x^{n} + y^{n} = z^{n}$ has no non-trivial solutions in
$\Q(\sqrt{2})$ provided $n \geq 4$. Their proof uses similar
ingredients to Wiles's work.

In \cite{DK}, Debarre and Klassen use Faltings's work on the rational
points on subvarieties of abelian varieties to prove that for $n \geq
3$ and $n \ne 6$, the equation $x^{n} + y^{n} = z^{n}$ has only
finitely many solutions $(x,y,z)$ where the variables belong to any
number field $K$ with $[K : \Q] \leq n-2$. Indeed, the work of Aigner shows
that when $n = 4$ the only non-trivial solution to $x^{n} + y^{n} = z^{n}$ with 
$x$, $y$ and $z$ in any quadratic field is
\[
  \left(\frac{1+\sqrt{-7}}{2}\right)^{4} + 
\left(\frac{1-\sqrt{-7}}{2}\right)^{4} = 1^{4},
\]
and when $n = 6$ or $n = 9$, there are no non-trivial solutions in quadratic 
fields.

We now turn to the problem of solutions to $x^{3} + y^{3} = z^{3}$ in
quadratic fields $\Q(\sqrt{d})$. For some choices of $d$ there are solutions,
such as
\[
  (18 + 17 \sqrt{2})^{3} + (18 - 17 \sqrt{2})^{3} = 42^{3}
\]
for $d = 2$, while for other choices (such as $d = 3$) there are no
non-trivial solutions. In 1913, Fueter \cite{Fueter} showed that if $d < 0$
and $d \equiv 2 \pmod{3}$, then there are no solutions if $3$ does not
divide the class number of $\Q(\sqrt{d})$. Fueter also proved in 
\cite{Fueter2} that there is a non-trivial solution to $x^{3} + y^{3}
= z^{3}$ in $\Q(\sqrt{d})$ if and only if there is one in
$\Q(\sqrt{-3d})$.

In 1915, Burnside \cite{Burnside} showed that every solution to $x^{3} +
y^{3} = z^{3}$ in a quadratic field takes the form
\begin{align*}
  x &= -3 + \sqrt{-3(1 + 4k^{3})},\\
  y &= -3 - \sqrt{-3(1 + 4k^{3})}\text{, and }\\ 
  z &= 6k
\end{align*}
up to scaling. Here $k$ is any rational number not equal to $0$ or $-1$.
This, however, does not answer the question of whether or not there are
solutions in $\Q(\sqrt{d})$ for given $d$ since it is not clear whether
\[
  dy^{2} = -3(1 + 4k^{3})
\]
has a solution with $k$ and $y$ both rational.

In a series of papers \cite{Aig2}, \cite{Aig1}, \cite{Aig4},
\cite{Aig3}, Aigner considered this problem (see \cite{Ribenboim},
Chapter XIII, Section 10 for a discussion in English). He showed that
there are no solutions in $\Q(\sqrt{-3d})$ if $d > 0$, $d \equiv 1
\pmod{3}$, and 3 does not divide the class number of
$\Q(\sqrt{-3d})$. He also developed general criteria to rule out the
existence of a solution. In particular, there are ``obstructing
integers'' $k$ with the property that there are no solutions in
$\Q(\sqrt{\pm d})$ if $d = kR$, where $R$ is a product of primes
congruent to 1 (mod 3) for which 2 is a cubic non-residue.

The goal of the present paper is to give a complete classification of
the fields $\Q(\sqrt{d})$ in which $x^{3} + y^{3} = z^{3}$ has a
solution. Our main result is the following.

\begin{thm}
\label{main}
Assume the Birch and Swinnerton-Dyer conjecture (see Section~\ref{BSD}
for the statement and background). If $d > 0$ is squarefree
with $\gcd(d,3) = 1$, then there is a non-trivial solution to
$x^{3} + y^{3} = z^{3}$ in $\Q(\sqrt{d})$ if and only if
\begin{align*}
  & \# \{ (x,y,z) \in \Z^{3} : x^{2} + y^{2} + 7z^{2} + xz = d \}\\
  &= \# \{ (x,y,z) \in \Z^{3} : x^{2} + 2y^{4} + 4z^{2} + xy + yz = d \}.
\end{align*}

If $d > 0$ is squarefree with $3 | d$, then there is a non-trivial solution
to $x^{3} + y^{3} = z^{3}$ in $\Q(\sqrt{d})$ if and only if
\begin{align*}
  & \# \{ (x,y,z) \in \Z^{3} : x^{2} + 3y^{2} + 27z^{2} = d/3 \}\\
  &= 
  \# \{ (x,y,z) \in \Z^{3} : 3x^{2} + 4y^{2} + 7z^{2} - 2yz = d/3 \}.
\end{align*}

Moreover, there are non-trivial solutions in $\Q(\sqrt{d})$ if and only if
there are non-trivial solutions in $\Q(\sqrt{-3d})$.
\end{thm}

\begin{rem}
Only one direction of our result is conditional on the
Birch and Swinnerton-Dyer conjecture. As mentioned in
Section~\ref{BSD}, it is known that if $E/\Q$ is an elliptic curve, 
$L(E,1) \ne 0$ implies that $E(\Q)$ is finite. As a consequence,
if the number of representations of $d$ (respectively $d/3$)
by the two different quadratic forms are different, then there are
no solutions in $\Q(\sqrt{d})$.
\end{rem}

Our method is similar to that used by Tunnell \cite{Tunnell} in his solution 
to the congruent number problem. The congruent number problem is to determine,
given a positive integer $n$, whether there is a right triangle with
rational side lengths and area $n$. It can be shown that $n$ is a congruent
number if and only if the elliptic curve $E_{n} : y^{2} = x^{3} - n^{2} x$
has positive rank. The Birch and Swinnerton-Dyer states that
$E_{n}$ has positive rank if and only if $L(E_{n},1) \ne 0$, and Waldspurger's
theorem (roughly speaking) states that
\[
  f(z) = \sum_{n=1}^{\infty} n^{1/4} \sqrt{L(E_{n},1)} q^{n}, \quad q = e^{2 \pi i z}
\]
is a weight $3/2$ modular form. Tunnell computes this modular form explicitly
as a difference of two weight $3/2$ theta series and proves that (in the case 
that $n$ is odd), $E_{n}$ is congruent if and only if $n$
has the same number of representations in the form $x^{2} + 4y^{2} + 8z^{2}$
with $z$ even as it does with $z$ odd. Tunnell's work was used in 
\cite{Trillion} to determine precisely which integers $n \leq 10^{12}$
are congruent (again assuming the Birch and Swinnerton-Dyer conjecture).

\begin{rem}
In \cite{Purkait}, Soma Purkait computes two (different) weight $3/2$ modular 
forms whose coefficients interpolate the central critical $L$-values of
twists of $x^{3} + y^{3} = z^{3}$ (see Proposition 8.7). Purkait expresses the 
first as a linear combiation of 7 theta series, but does not express the
second in terms of theta series. 
\end{rem}

An outline of the paper is as follows. In Section~\ref{BSD} we will 
discuss the Birch and Swinnerton-Dyer conjecture. In Section~\ref{background}
we will develop the necessary background. This will be used in
Section~\ref{proofs} to prove Theorem~\ref{main}.

\begin{ack}
  This work represents the master's thesis of the first author which
  was completed at Wake Forest University. The authors used Magma
  \cite{Magma} version 2.17 for computations in spaces of modular
  forms of integer and half-integer weights.
\end{ack}

\section{Elliptic Curves and the Birch and Swinnerton-Dyer 
Conjecture}\label{BSD}

The smooth, projective curve $C: x^{3} + y^{3} = z^{3}$ is an elliptic curve.
Specifically, if $X = \frac{12z}{y+x}$ and $Y = \frac{36(y-x)}{y+x}$, then
\[
  E : Y^{2} = X^{3} - 432.
\]
From Euler's proof of the $n = 3$ case of Fermat's last theorem, it follows
that the only rational points on $x^{3} + y^{3} = z^{3}$ are
$(1 : 0 : 1)$, $(0 : 1 : 1)$, and $(1 : -1 : 0)$. These correspond to the
three-torsion points $(12,-36)$, $(12,36)$, and the point at infinity
on $E$.

Suppose that $K = \Q(\sqrt{d})$ is a quadratic field and $\sigma : K
\to K$ is the automorphism given by $\sigma(a + b \sqrt{d}) = a - b
\sqrt{d}$ with $a, b \in \Q$. If $P = (x,y) \in E(K)$, define
$\sigma(P) = (\sigma(x), \sigma(y)) \in E(K)$. Then, $Q = P -
\sigma(P) \in E(K)$ and $\sigma(Q) = -Q$. Since the inverse of $(x,y)
\in E(K)$ is $(x,-y)$, it follows that $P - \sigma(P) = (a, b
\sqrt{d})$ for $a, b \in \Q$. Thus, $(a,b)$ is a rational point on the
quadratic twist $E_{d}$ of $E$, given by
\[
  E_{d} : dY^{2} = X^{3} - 432.
\]
\begin{lem}
\label{torslem}
The point $(a,b)$ on $E_{d}(\Q)$ is in the torsion subgroup of $E_{d}(\Q)$ if 
and only if the corresponding solution to $x^{3} + y^{3} = z^{3}$ is trivial.
\end{lem}
This lemma will be proven in Section~\ref{proofs}. Thus, there is a
non-trivial solution in $\Q(\sqrt{d})$ if and only if $E_{d}(\Q)$ has
positive rank.

If $E/\Q$ is an elliptic curve, let
\[
  L(E,s) = \sum_{n=1}^{\infty} \frac{a_{n}(E)}{n^{s}}
\]
be its $L$-function (see \cite{Silverman}, Appendix C, Section 16 for
the precise definition). It is known (see \cite{BCDT}) that
$L(E,s) = L(f,s)$ for some weight 2 modular form $f \in
S_{2}(\Gamma_{0}(N))$, where $N$ is the conductor of $E$. It follows
from this that $L(E,s)$ has an analytic continuation and functional
equation of the form
\[
  \Lambda(E,s) = (2 \pi)^{-s} N^{s/2} \Gamma(s) L(E,s)
\]
and $\Lambda(E,s) = w_{E} \Lambda(E,2-s)$, where
$w_{E} = \pm 1$ is the root number of $E$. Note that if
$w_{E} = -1$, then $L(E,1) = 0$.
The weak Birch and Swinnerton-Dyer conjecture predicts that
\[
  \ord_{s=1} L(E,s) = \rank(E(\Q)).
\]
The strong form predicts that
\[
  \lim_{s \to 1} \frac{L(E,s)}{(s-1)^{r}} = \frac{\Omega(E)
  R(E/\Q) \prod_{p} c_{p} \# \Sha(E/\Q)}{(\# E_{\tors})^{2}}.
\]
Here, $\Omega(E)$ is the real period of $E$ times
the number of connected components of $E(\R)$, $R(E/\Q)$ is the elliptic
regulator, the $c_{p}$ are the Tamagawa numbers, and $\Sha(E/\Q)$
is the Shafarevich-Tate group.

Much is known about the Birch and Swinnerton-Dyer in the case
when $\ord_{s=1} L(E,s)$ is $0$ or $1$. See for example 
\cite{CoatesWiles}, \cite{GrossZagier}, \cite{Kolyvagin}, and \cite{Rubin}.
The best known result currently is the following.

\begin{thm}[Gross-Zagier, Kolyvagin, et. al.]
Suppose that $E/\Q$ is an elliptic curve and $\ord_{s=1} L(E,s) = 0$ or $1$.
Then, $\ord_{s=1} L(E,s) = \rank(E(\Q))$.
\end{thm}
The work of Bump-Friedberg-Hoffstein \cite{BFH} or Murty-Murty \cite{MM}
is necessary to remove a condition imposed in the work of Gross-Zagier
and Kolyvagin.

\section{Preliminaries}
\label{background}

If $d$ is an integer, let $\chi_{d}$ denote the unique primitive
Dirichlet character with the property that
\[
  \chi_{d}(p) = \legen{d}{p}
\]
for all odd primes $p$. This character will be denoted by
$\chi_{d}(n) = \legen{d}{n}$, even when $n$ is not prime.

If $\lambda$ is a positive integer, let $M_{2\lambda}(\Gamma_{0}(N),\chi)$ 
denote the $\C$-vector space of modular forms of weight $2\lambda$ for 
$\Gamma_{0}(N)$ with character $\chi$, and $S_{2\lambda}(\Gamma_{0}(N), \chi)$
denote the subspace of cusp forms. Similarly, if $\lambda$ is a 
positive integer, let $M_{\lambda + \frac{1}{2}}(\Gamma_{0}(4N), \chi)$ denote the 
vector space of modular forms of weight $\lambda + \frac{1}{2}$ on
$\Gamma_{0}(4N)$ with character $\chi$ and $S_{\lambda + \frac{1}{2}}(\Gamma_{0}(4N), \chi)$ denote the subspace of cusp forms. We will frequently use
the following theorem of Sturm \cite{St} to prove that two modular forms
are equal.
\begin{thm}\label{sturm}
Suppose that $f(z) \in M_{r}(\Gamma_{0}(N), \chi)$ is a modular
form of integer or half-integer weight with 
$f(z) = \sum_{n=0}^{\infty} a(n) q^{n}$. If $a(n) = 0$ for 
$n \leq \frac{r}{12} [\SL_{2}(\Z) : \Gamma_{0}(N)]$, then $f(z) = 0$.
\end{thm}

We denote by $T_{p}$
the usual index $p$ Hecke operator on $M_{2 \lambda}(\Gamma_{0}(N), \chi)$,
and by $T_{p^{2}}$ the usual index $p^{2}$ Hecke operator on 
$M_{\lambda + 1/2}(\Gamma_{0}(4N), \chi)$.

Next, we recall the Shimura correspondence.
\begin{thm}[\cite{Shi}]
Suppose that $f(z) = \sum_{n=1}^{\infty} a(n) q^{n} \in
S_{\lambda + 1/2}(\Gamma_{0}(4N), \chi)$. For each squarefree integer $t$, 
let
\[
  \mathcal{S}_{t}(f(z)) = \sum_{n=1}^{\infty} 
  \left(\sum_{d | n} \chi(d) \legen{(-1)^{\lambda} t}{d}
  d^{\lambda - 1} a(t(n/d)^{2})\right) q^{n}.
\]
Then, $\mathcal{S}_{t}(f(z)) \in M_{2\lambda}(\Gamma_{0}(2N), \chi^{2})$. 
\end{thm}

One can show using the definition that if $p$ is a prime and
$p \nmid 4tN$, then
\[
  \mathcal{S}_{t}(f | T_{p^{2}}) = S_{t}(f) | T_{p},
\]
that is, the Shimura correspondence commutes with the Hecke action.

In \cite{Wald}, Waldspurger relates the Fourier coefficients of
a half-integer weight Hecke eigenform $f$ with the central critical $L$-values
of the twists of the corresponding integer weight modular form $g$
with the same Hecke eigenvalues. To state it, let $\Q_{p}$ be the
usual field of $p$-adic numbers. Also, if
\[
  F(z) = \sum_{n=1}^{\infty} A(n) q^{n},
\]
let $(F \otimes \chi)(z) = \sum_{n=1}^{\infty} A(n) \chi(n) q^{n}$.

\begin{thm}[\cite{Wald}, Corollaire 2, p. 379]\label{waldspurger}
Suppose that $f \in S_{\lambda + 1/2}(\Gamma_{0}(4N), \chi)$ is a half-integer
weight modular form and $f | T_{p^{2}} = \lambda(p) f$ for all $p \nmid 4N$.
Denote the Fourier expansion of $f(z)$ by 
\[
  f(z) = \sum_{n=1}^{\infty} a(n) q^{n}, \quad q = e^{2 \pi i z}.
\]
If $F(z) \in S_{2 \lambda}(\Gamma_{0}(2N), \chi^{2})$ is an integer weight modular 
form with $F(z) | T_{p} = \lambda(p) g$ for all $p \nmid 4N$ and $n_{1}$ and 
$n_{2}$ are two squarefree positive integers with 
$n_{1}/n_{2} \in \left(\Q_{p}^{\times}\right)^{2}$ for all $p | N$, then
\[
  a(n_{1})^{2} L(F \otimes \chi^{-1} \chi_{n_{2} \cdot (-1)^{\lambda}}, \lambda)
  \chi(n_{2}/n_{1}) n_{2}^{\lambda - 1/2}
  = a(n_{2})^{2} L(F \otimes \chi^{-1} \chi_{n_{1} \cdot (-1)^{\lambda}}, \lambda)
  n_{1}^{\lambda - 1/2}.
\]
\end{thm}

Our goal is to construct two modular forms 
$f_{1}(z) \in S_{3/2}(\Gamma_{0}(108))$ and $f_{2}(z) \in S_{3/2}(\Gamma_{0}(108),
\chi_{3})$ that  have the same Hecke eigenvalues as
\[
  F(z) = q \prod_{n=1}^{\infty} (1-q^{3n})^{2} (1-q^{9n})^{2} \in S_{2}(\Gamma_{0}(27)).
\]
This is the weight 2 modular form corresponding to $E_{1} : y^{2} = x^{3} - 432$.
As in \cite{Tunnell}, we will express $f_{1}$ and $f_{2}$ as linear combinations
of ternary theta functions. The next result recalls the modularity of
the theta series of positive-definite quadratic forms.

\begin{thm}[Theorem 10.9 of \cite{Iwaniec}]
\label{thetatrans}
Let $A$ be a $r \times r$ positive-definite symmetric matrix with
integer entries and even diagonal entries. Let $Q(\vec{x}) =
\frac{1}{2} \vec{x}^{T} A \vec{x}$, and let
\[
  \theta_{Q}(z) = \sum_{n=0}^{\infty} r_{Q}(n) q^{n}
\]
be the generating function for the number of representations of $n$ by $Q$. 
Then, 
\[
  \theta_{Q}(z) \in M_{r/2}(\Gamma_{0}(N), \chi_{\det(2A)}),
\]
where $N$ is the smallest positive integer so that
$NA^{-1}$ has integer entries and even diagonal entries.
\end{thm}

Finally, we require some facts about the root numbers of the curves
$E_{d}$. If $F(z) \in S_{2}(\Gamma_{0}(N))$ is the modular form corresponding 
to $E$, let $F(z) | W(N) = N^{-1} z^{-2} f\left(-\frac{1}{Nz}\right)$.
Then $F(z) | W(N) = -w_{E} F(z)$ 
(see for example Theorem 7.2 of \cite{Iwaniec}). Theorem 7.5 of
\cite{Iwaniec} states that if $\psi$ is a quadratic Dirichlet character with
conductor $r$ and $\gcd(r,N) = 1$, then $F \otimes \psi \in S_{2}(\Gamma_{0}(Nr^{2}))$ and
\[
  (F \otimes \psi) | W(Nr^{2}) = (\psi(N) \tau(\psi)^{2}/r) F | W(N)
\]
where $\tau(\psi) = \sum_{m=1}^{r} \psi(m) e^{2 \pi i m / r}$ is the usual
Gauss sum. 

Suppose $d$ is an integer so that $|d|$ is the conductor of $\chi_{d}$ and 
$F(z) \in S_{2}(\Gamma_{0}(27))$ is the modular form corresponding to $E_{1}$.
Then $F \otimes \chi_{d}$ is the modular form corresponding to $E_{d}$. 
Using the result from the previous paragraph and the equality
$\tau(\chi_{d})^{2} = |d| \chi_{d}(-1)$, we get
\[
  w_{E_{d}} = w_{E_{1}} \chi_{d}(27) \chi_{d}(-1) = \chi_{d}(-27).
\]
provided $\gcd(d,3) = 1$.

\section{Proofs}
\label{proofs}

In this section, we prove Lemma~\ref{torslem} and Theorem~\ref{main}.

Before we prove Lemma~\ref{torslem}, we will first need to determine the order of torsion
subgroup of $E_d(\mathbb{Q})$. First, note that $6\sqrt{2} \notin \mathbb{Z}$
and hence $E_d(\mathbb{Q})$ has no element of order two. Since there are no elements of 
order 2 in $E_d(\Q)_{\rm{tors}}$, then $2 \nmid |E_d(\Q)_{\rm{tors}}|$. We will now show $q \nmid E_d(\Q)_{\rm tors}$
 for primes $q > 3$.

If $p$ is prime with $p \equiv 2 \pmod{3}$, then we have that the map $x \rightarrow x^3 \in \F_{p}$
 is a bijection. Since this is a bijection, we have that
 $\sum_{x = 0}^{p-1} \legen{f(x)}{p} = 0$. Thus we have $\#E(\F_p) = p + 1$.
Suppose that $|E_d(\Q)_{\rm{tors}}| = N$ for $N$ odd.

If we suppose that a prime $q > 3$ divides $N$ then we can find an integer $x$ that is
relatively prime to $3q$ so that $x \equiv 2 \pmod{3}$ and $x \equiv 1 \pmod{q}$. By 
Dirichlet's Theorem, we have an infinite number of primes contained in the arithmetic progression
$3nq + x$ for $n \in \N$. If we take $p$ to be a sufficiently large prime in this progression, then the 
reduction of $E_d(\Q)_{\rm tors} \subseteq E(\F_{p})$  has order $N$. So, now we 
have $q | |E_d(\F_p)| = p + 1 \equiv x + 1 \equiv 2 \pmod{q}$. This is a contradiction. Hence
the only prime that divides $N$ is 3. We can follow a similar argument to show that 9 does not
divide $N$. This means that the torsion subgroup of $E_d(\Q)$ is either $\Z/3\Z$ or trivial.

Futhermore, if $E_d(\Q)$ contains a point of order 3 then the $x-$coordinate of the point
must be a root of the three-division polynomial $\phi_3(x) = 3x^4 - 12(432)d^3x$. The only real roots to $\phi_3(x)$ 
are $x = 0$ and $x = 12d$. For $x = 0$, then we have $y = \pm 108$ and $d = -3$. Finally
for $x = 12d$, then we find that $y = 1296d^3$ and that $d = 1$. Thus we conclude that
$E_d(\Q) \cong \Z/3\Z$ if and only if $d \in \{1, -3\}$. Finally, the torsion subgroup of $E_d(\Q)$
is trivial for $d \notin \{1, -3\}.$

\begin{proof}[Proof of Lemma 2]
$(\Rightarrow)$ Let $(x,y) \in E_d(\mathbb{Q})$ so that $(x,y)$ is not in 
$E_d(\mathbb{Q})_{\rm tors}$. By doing some arithmetic we get that 
$(x, y\sqrt{d}) \in E(K)$. In Section$~\ref{BSD}$, we defined a map
from $C(K) \rightarrow E(K)$. The inverse of this map sends 
$$(x,y\sqrt{d}) \rightarrow \left(\frac{\frac{12}{x} + \frac{y\sqrt{d}}{3x}}{2}, \frac{\frac{12}{x}-\frac{y\sqrt{d}}{3x}}{2}\right) \in C(K).$$
If we suppose that this is a trivial solution to $C$, then either the $x-$coordinate or $y-$coordinate
 is zero. Hence $y = \pm \frac{36\sqrt{d}}{d}$.

If $d = 1$, then we have $y = -36$ and $x = 12$. From Section 2, we know that $(12, -36)$ corresponds
to $(1:0:1)$ which is a trivial solution to $C$. Hence the point $(x,y)$ does not
satisfy the hypothesis for $d = 1$. Now for $d \neq 1$, we have $y \notin \mathbb{Q}$. This contradicts the hypothesis for $d \neq 1$.
Hence the solution we have is non-trivial.

$(\Leftarrow)$ Let $(x,y,z)$ be a non-trivial solution to $x^3 + y^3 = z^3$ in $K$.
Note that for $d = 1 \text{ or } -3$ Euler showed that there are only trivial solutions
and thus this direction is vacuously true for these two cases.

For $d \neq 1 \text{ and } -3$, from Section$~\ref{BSD}$ we showed that $(x,y,z) \rightarrow (X,Y) = P \in E(K)$.
 Also from section 2, if $P - \sigma(P) = (a, b\sqrt{d})$ then $(a, b) \in E_d(\mathbb{Q})$. Since $d \neq 1$ 
and $-3$, then the torsion subgroup of $E_d(\mathbb{Q})$ is trivial. Thus $(a, b) \notin E_d(\mathbb{Q})_{\rm{tors}}$. 
\end{proof}

Recall from Section$~\ref{background}$ that the elliptic curve $E$ corresponds to the modular form
 $F(z) = q\prod^\infty_{n = 1} (1 - q^{3n})^2(1-q^{9n})^2 \in S_2(\Gamma_0(27)).$

\begin{rem}
For convenience we will think of $F(z)$ as a Fourier series with coefficients $\lambda(n)$ for $n \in \N$. Note that if $\lambda(n) 
\neq 0$ then $n \equiv 1 \pmod{3}$. So we can write $\lambda(n) = \lambda(n) (\frac{n}{3})$ for $n \in \N$. Hence $F\otimes \chi_{-3d} =
 F \otimes \chi_d$. We can now conclude that $L(E_d, 1) = L(E_{-3d},1)$.
\end{rem}

\begin{proof}[Proof of Theorem 1]
To begin we will examine the case for $d < 0$ so that $d \equiv 2 \pmod{3}$. Note that ${\rm dim}~S_{3/2}(\Gamma_0(108),\chi_1) = 5$.
 Moreover, we have the following basis of $S_{3/2}(\Gamma_0(108),\chi_1))$:
\begin{align*}
    g_1(z) &= q - q^{10} - q^{16} - q^{19} - q^{22} + 2q^{28} + \cdots,\\
    g_2(z) &= q^{2} - q^{5} + q^{8} - q^{11} + q^{14}- 2q^{17} - q^{20} + \cdots,\\
    g_3(z) &= q^{3} - 2q^{12} + \ldots,\\
    g_4(z) &= q^{4} - q^{10} + q^{13} - q^{16} - q^{19} - q^{22} - q^{25} + q^{28} + \cdots \text{, and}\\
    g_5(z) &= q^{7} - q^{10} + q^{13} - q^{16} - q^{22} - q^{25} + \cdots.
\end{align*}

By Theorem$~\ref{sturm}$ we have:
\begin{align*}
\mathcal{S}_1(g_1(z) + g_4(z)) &= F(z) + F(z)|V(2),\\
\mathcal{S}_2(g_1(z) + g_4(z)) &= 0,\\
\mathcal{S}_3(g_1(z) + g_4(z)) &= 0,\\
\mathcal{S}_1(g_1(z) + g_5(z)) &= F(z), \\
\mathcal{S}_2(g_1(z) + g_5(z)) &= 0 \text{, and}\\
\mathcal{S}_3(g_1(z) + g_5(z)) &= 0.
\end{align*}

Since we took $t = 1, 2,$ and $3$, then from Section$~\ref{BSD}$ we have $\mathcal{S}_t((g_1(z) + g_4(z))|T_{p^2}) = \mathcal{S}_t((g_1(z) + g_4(z))|T(p)$ for all primes $p > 3$. Since $F(z)$ and $F(z)|V(2)$ are both Hecke eigenforms, then $F(z)|T(p) = \lambda(p)F(z)$ and $F(z)|V(2)|T(p) = \lambda(p)F(z)|V(2)$ for primes $p > 3$. Also, since 1, 2 and 3 divide $4N$, then we have $(g_1(z) + g_4(z)|T_{p^2} - \lambda(p)(g_1(z) + g_4(z)$ is in: $\rm{ker}(\mathcal{S}_1)$, $\rm{ker}(\mathcal{S}_2)$, and $\rm{ker}(\mathcal{S}_3)$. Furthermore since $\rm{ker}(\mathcal{S}_1) \cap \rm{ker}(\mathcal{S}_2) \cap \rm{ker}(\mathcal{S}_3) = 0$, then $g_1(z) + g_4(z)$ is a Hecke eigenform. The case of $g_1(z) + g_5(z)$ is similar.

We will now take the quadratic forms $Q_1(x,y,z) = x^2 + 3y^2 + 27z^2$ and $Q_2(x,y,z) = 3x^2 + 4y^2 - 2yz + 7z^2$.
We have their theta-series $\theta_{Q_1}, \theta_{Q_2} \in M_{3/2}(\Gamma_0(108), \chi_1)$. Also by Theorem$~\ref{sturm}$, we have
 $$\theta_{Q_1}(z) - \theta_{Q_2}(z) = -2(g_1(z) + g_5(z)) + 4(g_1(z) + g_4(z)).$$
Furthermore, since $g_1(z) + g_4(z)$ and $g_1(z) + g_5(z)$ are both Hecke eigenforms with the same eigenvalues 
then $\theta_{Q_1}(z) - \theta_{Q_2}(z)$ is a Hecke eigenform as well.

Let $a(n)$ denote the $n-$th cofficient of $\theta_{Q_1}(z) - \theta_{Q_2}(z)$. By Theorem$~\ref{waldspurger}$, we have
$$L(E_{-n_2},1) = \sqrt{\frac{n_2}{n_1}}\left(\frac{a(n_2)}{a(n_1)}\right)^2 L(E_{-n_1},1)$$
for $n_1$ and $n_2$ squarefree with $\legen{n_1/n_2}{p} = 1$ for $p = 3$ and $n_1/n_2 \equiv 1 \pmod{8}$. If
 we take $n_2 \equiv 1 \pmod{3}$, then the table below covers all possible cases.

\begin{center}
\begin{tabular}{cccccc}
$n_2$ & $n_1$ & $a(n_1)$ & $L(E_{-n_1},1)$\\
\hline\\
$n_2 \equiv 1 \pmod{24}$    &1     &  $2 $   & $1.52995 \ldots$\\
$n_2 \equiv 34 \pmod{48}$  &34   &  $4 $   & $1.04953 \ldots$\\
$n_2 \equiv 19 \pmod{24}$  &19   & $-6$    & $0.70199 \ldots$\\
$n_2 \equiv 13 \pmod{24}$  &13   & $ 2$    & $0.42434 \ldots$\\
$n_2 \equiv 22 \pmod{48}$  &22   & $-4$   & $1.30474 \ldots$\\
$n_2 \equiv 7 \pmod{24}$    &7     & $-2$   & $1.15653 \ldots$\\ 
$n_2 \equiv 10 \pmod{36}$  &10   & $-4$  & $1.93525 \ldots$\\
$n_2 \equiv 46 \pmod{48}$  &46   & $4$   & $0.90231 \ldots$
\end{tabular}
\end{center}

Thus we have $d < 0$ with $d \equiv 2 \pmod{3}$ and $L(E_d, 1) = 0$ if and only if $a(-d) = 0$. Since $L(E_{3n_2},1) = L(E_{-n_2}, 1)$ then we have $d > 0$ so that $d \equiv 3 \pmod{9}$ and $L(E_d, 1) = 0$ if and only if $a(d/3) = 0$.

We will now examine the case $d < 0$ so that $3|d$ and $d \equiv 6 \pmod{9}$. Note that ${\rm{dim}}~ S_{3/2}(\Gamma_0(108), \chi_3) = 5$, and we have the basis:
\begin{align*}
   h_1(z) &=  q - 2q^{13} - q^{25} - 2q^{28} + \cdots, \\
   h_2(z) &= q^{2} + q^{5} - q^{8} - q^{11} - q^{14} - q^{20} - 2q^{23} + 2q^{26} + \cdots, \\
   h_3(z) &= q^{4} - q^{13} - 2q^{16} + 2q^{25} - q^{28} + \cdots, \\
   h_4(z) &=  q^{7} - q^{13} - q^{19} + \cdots \text{, and}\\
   h_5(z) &= q^{10} - q^{16} - q^{19} - q^{22} + q^{25} + \cdots.
\end{align*}

By Theorem$~\ref{sturm}$, we have
\begin{align*}
\mathcal{S}_1(h_1(z) - h_4(z) + 2h_5(z)) &= F(z), \\
\mathcal{S}_2(h_1(z) - h_4(z) + 2h_5(z)) &= 0, \\
\mathcal{S}_3(h_1(z) - h_4(z) + 2h_5(z)) &= 0, \\
\mathcal{S}_1(h_1(z) -4h_3(z) - 5h_4(z) + 10h_5(z)) &= F(z) + 4F(z)|V(2),\\
\mathcal{S}_2(h_1(z) -4h_3(z) - 5h_4(z) + 10h_5(z)) &= 0 \text{, and}\\ 
\mathcal{S}_3(h_1(z) -4h_3(z) - 5h_4(z) + 10h_5(z)) &= 0.
\end{align*}

From a similar argument as the previous case, we get that $h_1(z) - h_4(z) + 2h_5(z)$ and  $h_1(z) - 4h_3(z) - 5h_4(z) + 10h_5(z)$ are Hecke eigenforms for $T(p^2)$ for primes $p > 3$. We will now take the quadratic forms $Q_3(x,y,z) = x^2 + y^2 + 7z^2 + xz$ and $Q_4(x,y,z) = x^2 + 2y^2 + 4z^2 + xy + yz$. 

We will denote the theta series corresponding to $Q_3$ and $Q_4$ by $\theta_{Q_3}$ and $\theta_{Q_4}$, respectively. Note that
 $\theta_{Q_3}, \theta_{Q_4} \in  M_{3/2}(\Gamma_0(108), \chi_3)$. 
By Theorem$~\ref{sturm}$, $\theta_{Q_3} - \theta_{Q_4} = 2h_1(z) - 4h_3(z) - 6h_4(z) + 12h_5(z)$. Since  $h_1(z) - h_4(z) + 2h_5(z)$ and  $h_1(z) - 4h_3(z) - 5h_4(z) + 10h_5(z)$ have the same eigenvalues, $\theta_{Q_3} - \theta_{Q_4}$ is a Hecke eigenform. Let $b(n)$
 denote the $n-$th coefficient of $\theta_{Q_3} - \theta_{Q_4}$.

Hence by Theorem$~\ref{waldspurger}$, we have $$L(E_{-3n_2}, 1) = \sqrt{\frac{n_2}{n_1}}\left(\frac{b(n_2)}{b(n_1)}\right)^2 L(E_{-3n_1}, 1).$$

\begin{center}
\begin{tabular}{cccccc}
$n_2$ & $n_1$ & $b(n_1)$ & $L(E_{-3n_1},1)$\\
\hline\\
$n_2 \equiv 1 \pmod{24}$    &1       &$2$       & $0.58887 \ldots$\\
$n_2 \equiv 34 \pmod{48}$  &34     &$12$     & $1.81785 \ldots$\\
$n_2 \equiv 19 \pmod{24}$  &19     &$-6$      & $0.60794 \ldots$\\
$n_2 \equiv 13 \pmod{24}$  &13     &$6$       & $1.46993 \ldots$\\
$n_2 \equiv 22 \pmod{48}$  &22     &$-12$    & $2.25989 \ldots$\\
$n_2 \equiv 7 \pmod{24}$    &7       &$-6$      & $1.00159 \ldots$\\ 
$n_2 \equiv 10 \pmod{36}$  &10     &$12$     & $3.35196 \ldots$\\  
$n_2 \equiv 46 \pmod{48}$  &46     &$-12$    & $1.56286 \ldots$
\end{tabular}
\end{center}

Therefore if $d < 0$ and $d \equiv 6 \pmod{9}$, $L(E_{d}, 1) = 0$ if and only if $b(-d/3) = 0$. Furthermore by the remark, we have $$L(E_{n_2}, 1) = \sqrt{\frac{n_2}{n_1}}\left(\frac{b(n_2)}{b(n_1)}\right)^2 L(E_{-3n_1}, 1)$$ for $n_2 \equiv 1 \pmod{3}$. Thus for squarefree $d > 0$ so that $d \equiv 1 \pmod{3}$, $L(E_d, 1) = 0$ if and only if $b(d) = 0$.

There are two pairs of cases that Theorem$~\ref{waldspurger}$ does not handle: $d > 0$ with $d \equiv 6 \pmod{9}$ and $d < 0$
 with $d \equiv 1 \pmod{3}$, and $d > 0$ with $d \equiv 2 \pmod{3}$ and $d < 0$ with $d \equiv 3 \pmod{9}$. 

We will now consider $d > 0$ with $d \equiv 2 \pmod{3}$. To handle this case, we will show that the root number of $E_d$ is $-1$.
Recall from the end of Section$~\ref{background}$ that $w_{E_d} = \chi_d(-27)$. Hence $w_{E_d} = \chi_d(-27) = \legen{d}{-27} = -1$. 
 Furthermore, by the remark we have $w_{E_d} = -1$ for $d < 0$ so that $d \equiv 3 \pmod{9}$.

We now want to show that there are no non-trivial solutions for $d > 0$ with $d \equiv 6 \pmod{9}$. To do this, we will show that $x^2 +
 3y^2 + 27z^2 = 3x^2 + 4y^2 - 2yz + 7y^2 \neq d$. Since $d > 0$ then $-d/3 \equiv 1 \pmod{3}$. This means that $x^2 + 3y^2 + 27z^2
 \equiv 0$ or $1 \pmod{3}$. We also have that $3x^2 + 4y^2 - 2yz + 7z^2 \equiv (y + 2z)^2 \equiv 0$ or $1 \pmod{3}$. Hence $r_{Q_1}(d) = r_{Q_2}(d) = 0$.

Let $\psi$ be the non-trivial Dirichlet character with modulus $3$. Note that $(\theta_{Q_3}(z) - \theta_{Q_4})\otimes \psi 
\in M_{3/2}(\Gamma_0(108*3^2), \chi_3 \psi^2)$ by Proposition 3.12 in \cite{Ono}.
By Theorem$~\ref{sturm}$, $(\theta_{Q_3}(z) - \theta_{Q_4}) \otimes \psi = \theta_{Q_3}(z) - \theta_{Q_4}(z)$. So $b(n)\psi(n) = b(n)$ for all $n \geq 1$. If $d \equiv 2 \pmod{3}$, we have that $\psi(d) = -1$. Thus $b(d) = -b(d)$. Therefore $b(d) = 0$. Thus $r_{Q_3}(d) = r_{Q_4}(d)$ for $d \equiv 2 \pmod{3}$.

Hence we have shown that by checking the number of solutions of the pair of equations 
 $x^{2} + y^{2} + 7z^{2} + xz$ and $x^{2} + 2y^{2} + 4z^{2} + xy + yz$, and
 $x^{2} + 3y^{2} + 27z^{2}$ and  $3x^{2} + 4y^{2} + 7z^{2} - 2yz  $
 is sufficient to determine when there are non-trivial solutions to $x^3 + y^3 = z^3$ in $\Q(\sqrt{d})$.
\end{proof}

\bibliographystyle{plain}
\bibliography{refs}

\end{document}